\input amstex
\documentstyle{amsppt}
\magnification=\magstep1
 \hsize 13cm \vsize 18.35cm \pageno=1
\loadbold \loadmsam
    \loadmsbm
    \UseAMSsymbols
\topmatter
\NoRunningHeads
\title  The modified $q$-Euler numbers and polynomials
\endtitle
\author
 Taekyun Kim
\endauthor
\abstract In the recent  paper (see [6]) we  defined a set of
numbers inductively by
$$
E_{0,q} =1, \ \ q(qE +1)^n + E_{n,q} = \cases
[2]_q , & \text{  if  $n = 0$}\\
0 ,&  \text{  if  $n \neq 0$,}
\endcases
$$
with the usual convention of replacing $E^n$ by $E_{n,q}$. These
numbers $E_{k,q}$ are called ``the $q$-Euler numbers" which are
reduced to $E_k$ when $q=1$. In this paper we construct the
modified $q$-Euler numbers  $\Cal E_{k,q}$
$$
\Cal E_{0,q} =\frac{[2]_q}{2} , \ \ (q \Cal E  +1)^n + \Cal
E_{n,q} = \cases
[2]_q , & \text{  if  $n = 0$,}\\
0 ,&  \text{  if  $n \neq 0$,}
\endcases
$$
with the usual convention of replacing $\Cal E^i$ by $\Cal
E_{i,q}$. Finally we give some interesting identities related to
these $q$-Euler numbers $\Cal E_{i,q}$.
\endabstract
\endtopmatter

\document

{\bf\centerline {\S 1. Introduction}}

 \vskip 20pt

Let $p$ be a fixed odd prime. Throughout this paper
 $\Bbb Z_p$, $\Bbb Q_p$,  $\Bbb C$ and $\Bbb C_p$ will
 respectively,
denote the ring of $p-$adic rational integers, the field of
$p$-adic rational numbers, the complex number field and the
completion of the algebraic closure of $\Bbb Q_p$. Let $\nu_p$ be
the normalized exponential valuation of $\Bbb C_p$ with $|p|_p
=p^{-\nu_p (p)} =\frac{1}{p}$. When one talks of $q$-extension,
$q$ is variously considered as an indeterminate, a complex
$q\in\Bbb C$, or a $p$-adic number $q\in\Bbb C_p$. If $q\in\Bbb
C$, one normally assumes $|q|<1$. If $q\in\Bbb C_p$, then we
assume $|q-1|_p < p^{-1/p-1}$, so that $q^x =\exp (x \log q)$ for
$|x|_p \leq 1$. The ordinary Euler numbers are defined by  the
generating function as follows:
$$
F(t) = \frac{2}{e^t +1} =\sum_{n=0}^\infty E_n \frac{t^n}{n!}, \
\text{ cf. [6].}
$$
From this equation, we derive the following relation:

$$
E_{0} =1, \ \ (E +1)^n + E_{n} = \cases
2 , & \text{  if  $n = 0$,}\\
0 ,&  \text{  if  $n \neq 0$,}
\endcases
$$
where we use the technique method notation by replacing $E^n$ by
$E_{n,q}$ ($n\geq 0$), symbolically. In the recent(see[6,8]), we
defined  ``the $q$-Euler numbers" as
$$\eqalignno{ &
E_{0,q} =1, \ \ q(qE +1)^n + E_{n,q} = \cases
[2]_q , & \text{  if  $n = 0$,}\\
0 ,&  \text{  if  $n \neq 0$,}
\endcases &(1)}
$$
with the usual convention of replacing $E^n$ by $E_{n,q}$. These
numbers are reduced to $E_k$ when $q=1$. From (1), we also derive
$$
E_{n,q} =\frac{[2]_q}{(1-q)^n} \sum_{l=0}^n \binom{n}{l}
\frac{(-1)^l}{1+q^{l+1}}, \ \ \text{( see [6]),}
$$
where $\binom{n}{l} =\frac{n(n-1)\cdots(n-l+1)}{l!}$. The
$q$-extension of $n\in\Bbb N$ is defined by
$$
[n]_q = \dfrac{1-q^n}{1-q} =1+q +q^2 +\cdots +q^{n-1},
$$
and
$$
[n]_{-q} = \dfrac{1-(-q)^n}{1+q} =1-q +q^2 -\cdots +(-q)^{n-1},\ \
\text{ cf. [4,5,7,9].}
$$
In [1,2], Carlitz defined a set of numbers $\xi_k =\xi_k (q)$
inductively by
$$\eqalignno{ &
\xi_{0} =1, \ \ (q\xi +1)^k - \xi_k = \cases
1 , & \text{  if  $k = 1$,}\\
0 ,&  \text{  if  $k >1$,}
\endcases &(2)}
$$
with the usual convention of replacing $\xi^i$ by $\xi_i$. These
numbers are $q$-extension of ordinary Bernoulli numbers $B_k$, but
they do not remain finite when $q=1$. So, he modified the Eq.(2)
as follows:
$$\eqalignno{ &
\beta_{0} =1, \ \ q(q\beta +1)^k - \beta_k = \cases
1 , & \text{  if  $k = 1$,}\\
0 ,&  \text{  if  $k >1$.}
\endcases &(3)}
$$
These numbers $\beta_k =\beta_k (q)$ are called ``the
$q$-Bernoulli numbers", which are reduced to $B_k$ when $q=1$,
see [1,2].  Some properties of $\beta_k$ were investigated by many
authors (see [1,2,3,4,9]). In [3,10], the definition of modified
$q$-Bernoulli numbers $B_{k,q}$ are introduced by
$$\eqalignno{ &
B_{0,q} = \frac{q-1}{\log q} , \ \ (qB +1)^k - B_{k,q} = \cases
1 , & \text{  if  $k = 1$,}\\
0 ,&  \text{  if  $k >1$,}
\endcases &(3)}
$$
with the usual convention of replacing $B^i$ by $B_{i,q}$. For a
fixed positive integer $d$ with $(p,d)=1$, set
$$\split
& X=X_d = \lim_{\overleftarrow{N} } \Bbb Z/ dp^N \Bbb Z ,\cr & \
X_1 = \Bbb Z_p , \cr  & X^\ast = \underset {{0<a<d p}\atop
{(a,p)=1}}\to {\cup} (a+ dp \Bbb Z_p ), \cr & a+d p^N \Bbb Z_p =\{
x\in X | x \equiv a \pmod{dp^N}\},\ \text{cf. [3,4,5,6,7,8],}
\endsplit$$ where $a\in \Bbb Z$ satisfies the condition  $0\leq a < d p^N$.

We say that $f$ is uniformly differentiable function at a point $a
\in\Bbb Z_p$, and write $f\in UD(\Bbb Z_p )$, if the difference
quotient $ F_f (x,y) = \dfrac{f(x) -f(y)}{x-y} $ have a limit
$f^\prime (a)$ as $(x,y) \to (a,a)$. For
 $f\in UD(\Bbb Z_p )$, an invariant $p$-adic  $q$-integral was defined by
$$
I_q (f) =\int_{\Bbb Z_p }f(x) d\mu_q (x) = \lim_{N\to \infty}
\dfrac{1}{[p^N ]_q} \sum_{x=0}^{p^N -1} f(x) q^x ,\quad \text{see
}[3,4].
$$

From this we can derive
$$\eqalignno{ &
q I_q (f_1 ) = I_q (f) +(q-1) f(0) + \frac{q-1}{\log q}f^\prime
(0), \ \ \text{ see [3],} &(5)}
$$
where $f_1 (x) =f(x+1)$, $f^\prime (0) = \frac{df(0)}{dx}$.

In the sense of fermionic, let us define
$$\eqalignno{ &
q I_{-q } (f_1 ) = \lim_{q\to -q} I_q (f) = \int_{\Bbb Z_p }f(x)
d\mu_{-q} (x), \ \ \text{ see [6].} &(6)}
$$

Thus, we have the following integral relation:
$$
q I_{-q} (f_1 ) + I_{-q } (f) =[2]_q f(0),
$$
where $f_1 (x) =f(x+1)$. Let $ I_{-1} (f) =\lim_{q\to 1} I_{-q}
(f). $ Then we see that
$$
\int_{\Bbb Z_p} e^{tx} d\mu_{-1} (x) = \frac{2}{e^t +1}
=\sum_{n=0}^\infty E_n \frac{t^n}{ n!}, \ \ \text{ see [6,8].}
$$
In the present paper we give a new construction of $q$-Euler
numbers which can be uniquely determined by

$$
\Cal E_{0,q} =\frac{[2]_q}{2}, \ \ (q \Cal E  +1)^n + \Cal E_{n,q}
= \cases
[2]_q , & \text{  if  $n = 0$,}\\
0 ,&  \text{  if  $n \neq 0$,}
\endcases
$$
with the usual convention of replacing $\Cal E^n$ by $\Cal
E_{n,q}$. These $q$-Euler numbers are corresponding to
$q$-Bernoulli numbers $B_{k,q}$. Finally we shall consider
$q$-zeta function which interpolates $\Cal E_{k,q}$ at negative
integers. As an application of these numbers $\Cal E_{k,q}$, we
will investigate some interesting alternating sums of powers of
consecutive $q$-integers.

{\bf\centerline {\S 2. A note on $q$-Bernoulli and Euler numbers}}
{\bf\centerline { associated with $p$-adic $q$-integrals on $\Bbb
Z_p$}}

In [4],  it was known that

$$\eqalignno{ &  \int_{\Bbb Z_p} e^{[x]_q t} d\mu_{q} (x)
=\sum_{n=0}^\infty \beta_n \frac{t^n}{n!},&(7)}
$$
where $\beta_n$ are Carlitz's $q$-Bernoulli numbers. By (5) and
(7), we see that
$$
(q-1) +t =qI_q (e^{(1+q[x]_q)t}) -I_q (e^{[x]_q t}), \ \text{ cf.
[3]. }
$$
Thus, we have
$$\eqalignno{ &
\beta_{0} =1, \ \ q(q\beta +1)^n - \beta_n = \cases
1 , & \text{  if  $n = 1$}\\
0 ,&  \text{  if  $n >1$.}
\endcases &}$$

\proclaim{Lemma 1} For $n\in \Bbb N$, we have
$$\eqalignno{ & q^n I_q (f_n ) = I_q (f) +(q-1)\sum_{l=0}^{n-1}
q^l f(l) + \frac{q-1}{\log q} \sum_{l=0}^{n-1} q^l f^\prime (l),
&(8)}$$
 where $f_n (x) = f(x+n)$.
\endproclaim

\demo{Proof} By Eq.(5) and  induction, Lemma 1 can be easily
proved.
\enddemo

It was known that

$$\eqalignno{ & \int_{\Bbb Z_p} [x+y]_q^n d\mu_q (x) =\beta_n (x)
,\quad \text{ see [3,4],} &(9)}$$
 where $\beta_n (x)$ are  Carlitz's $q$-Bernoulli polynomials. Let
 $n\in\Bbb N$, $k\in\Bbb Z_+$. Then, by (8), we have

$$\eqalignno{ & q^n \int_{\Bbb Z_p} [x+n]_q^k d\mu_q (x)  = \int_{\Bbb Z_p} [x]_q^k d\mu_q (x)
+
 (q-1)\sum_{l=0}^{n-1}
q^l [l]_q^k + k \sum_{l=0}^{n-1} q^{2l} [l]_q^{k-1}. &(10)}$$
 By (7), (8) and (10), we obtain the following:

 \proclaim{Proposition 2}
For $n\in\Bbb N$, $k\in\Bbb Z_+$, we have
$$
q^n \beta_k (n) -\beta_k = (q-1)\sum_{l=0}^{n-1} q^l [l]_q^k + k
\sum_{l=0}^{n-1} q^{2l} [l]_q^{k-1} .
$$
  \endproclaim

\remark{Remark}  If we take $n=1$ in Proposition 2,  then we have
$$
q(q\beta +1)^k -\beta_k =\cases
1 , & \text{  if  $k = 1$,}\\
0 ,&  \text{  if  $k>1$.}
\endcases
$$
\endremark

In [10], it was also known that the modified $q$-Bernoulli numbers
and polynomials can be represented by $p$-adic $q$-integral as
follows:

$$\eqalignno{ & \int_{\Bbb Z_p}q^{-x} [x]_q^n d\mu_q (x)  = B_{n,q}, \ \text{ and }
\  \int_{\Bbb Z_p} q^{-y} [y+x]_q^n d\mu_{-q} (y) =B_{n,q} (x).
&(11)}$$

From the definition of $p$-adic $q$-integral, we easily derive
$$\eqalignno{ & I_q (q^{-x}f_1) =I_q (q^{-x} f) + \frac{q-1}{\log q}f^\prime
(0). &(12)}$$

By (12), we obtain the following lemma:

\proclaim{Lemma 3}  For $n\in\Bbb N$, we have
$$\eqalignno{ & I_q (q^{-x}f_n ) =I_q (q^{-x} f) + \frac{q-1}{\log q}\sum_{l=0}^{n-1}
f^\prime(l) , &(13)}
$$  where $f_n (x) =f(x+n).$ That is,
$$
\int_{\Bbb Z_p} q^{-x} f(x+n) d\mu_q (x) =\int_{\Bbb Z_p} q^{-x}
f(x) d\mu_q (x) + \frac{q-1}{\log q} \sum_{l=0}^{n-1} f^\prime
(l).
$$
\endproclaim

From (11) and  (13), we note that
$$\eqalignno{ & B_{k,q} (n) -  B_{k,q} =k\sum_{l=0}^{n-1} q^l
[l]_q^{k-1}, \ \ \text{ cf. [5,7],} &(14)}
$$
where $n\in\Bbb N$, $k\in\Bbb Z_+$. By the definition of $I_{-q}
(f)$, we show that
$$\eqalignno{ &  qI_{-q} (f_1 ) + I_{-q} (f) =[2]_q f(0). &(15)}
$$

From (15) and induction, we derive the following integral
equation:

$$\eqalignno{ &  q^n I_{-q} (f_n ) +(-1)^{n-1} I_{-q} (f) =[2]_q \sum_{l=0}^{n-1} (-1)^{n-1-l} q^l f(l), &(16)}
$$
where $n\in\Bbb N$, $f_n (x) =f(x+n )$. When $n$ is an odd
positive integer, we have
$$\eqalignno{ &  q^n I_{-q} (f_n ) + I_{-q} (f) =[2]_q \sum_{l=0}^{n-1} (-1)^{l} q^l f(l). &(17)}
$$
If $n$ is an even natural number, then we see that
$$\eqalignno{ &  q^n I_{-q} (f_n ) - I_{-q} (f) =[2]_q \sum_{l=0}^{n-1} (-1)^{l-1} q^l f(l). &(18)}
$$
By (17) and (18), we obtain the following lemma:

\proclaim{Lemma 4} Let $n$ be an odd positive integer. Then
$$[2]_q\sum_{l=0}^{n-1} q^l [l]_q^m =q^n E_{m,q} (n) +E_{m,q}.
$$
If $n$(=even)$\in \Bbb N$, then we have
$$
q^n E_{m,q} (n) -E_{m,q}  =[2]_q
 \sum_{l=0}^{n-1} (-1)^{l-1} q^l [l]_q^m .$$
\endproclaim

Let us consider the modified $q$-Euler numbers and polynomials.
For any non-negative integer $n$, the modified $q$-Euler numbers
$\Cal E_{n,q}$ are defined by

$$\Cal E_{n,q} =\int_{\Bbb Z_p} q^{-x} [x]_q^n d\mu_{-q} (x)
=[2]_q (\frac{1}{1-q} )^n \sum_{l=0}^n \binom{n}{l} (-1)^l
\frac{1}{1+q^l }.
$$
By using $p$-adic $q$-integral on $\Bbb Z_p$, we can also consider
the modified $q$-Euler polynomials as follows:

$$\Cal E_{n,q} (x)  =\int_{\Bbb Z_p} q^{-t} [x+t]_q^n d\mu_{-q} (t)
=[2]_q (\frac{1}{1-q} )^n \sum_{l=0}^n \binom{n}{l} (-1)^l
\frac{q^{xl}}{1+q^l } .
$$

From (6) and (15), we derive the following $p$-adic $q$-integral
relation:

$$\eqalignno{ &   I_{-q} (q^{-x}f_1 ) + I_{-q} (q^{-x} f) =[2]_q f(0).
&(19)}$$

Thus, we obtain the following proposition:

\proclaim{Proposition 5} For $n\in\Bbb N$, we have
$$\eqalignno{ &   I_{-q} (q^{-x}f_n ) +(-1)^{n-1} I_{-q} (q^{-x}f) =[2]_q \sum_{l=0}^{n-1} (-1)^{n-l-1}  f(l). &(20)}
$$
\endproclaim

If we take $f(x) =e^{[x]_q t}$ in Eq.(19), then we have
$$\eqalignno{  [2]_q &=
 \int_{\Bbb Z_p}   q^{-x}e^{[x+1]_q t}d\mu_{-q}(x)+ \int_{\Bbb Z_p}   q^{-x}e^{[x]_q
 t}d\mu_{-q}(x)\cr
 &=\sum_{n=0}^\infty  ( \sum_{l=0}^n \binom{n}{l}q^l \Cal E_{l,q}  + \Cal E_{l,q}
 )\frac{t^n}{n!} &(21)\cr
 &=\sum_{n=0}^\infty  ( (q\Cal E +1)^n   + \Cal E_{n,q}
 )\frac{t^n}{n!},}
$$
with the usual convention of replacing $\Cal E^n$ by $\Cal
E_{n,q}$.

Therefore we obtain the following theorem:

\proclaim{Theorem 6} Let  $n\in\Bbb Z_+$. Then
$$
(q\Cal E +1)^n + \Cal E_{n,q} = \cases
[2]_q , & \text{  if  $n = 0$,}\\
0 ,&  \text{  if  $n \neq 0$,}
\endcases
$$
with the usual convention of replacing $\Cal E^i$ by $\Cal
E_{i,q}$.
\endproclaim

Note that $\lim_{q\to 1} \Cal E_{n,q} =E_n$, where $E_n$ are
ordinary Euler numbers. From (19) and (20), we can derive the
following theorem:

\proclaim{Theorem 7} Let  $k$(=even)$\in\Bbb N$, and let $n\in\Bbb
Z_+$.  Then we have
$$ \Cal E_{n,q} - \Cal E_{n,q} (k) =[2]_q \sum_{l=0}^{k-1} (-1)^l
[l]_q^n.
$$
If  $k$(=odd)$\in\Bbb N$ and  $n\in\Bbb Z_+$, then we see that
$$ \Cal E_{n,q} + \Cal E_{n,q} (k) =[2]_q \sum_{l=0}^{k-1} (-1)^l
[l]_q^n.
$$
\endproclaim

Let $\chi$ be the Dirichlet's character with conductor
$d$(=odd)$\in \Bbb N$. Then we consider the modified generalized
$q$-Euler numbers attached to $\chi$ as follows:

$$
\Cal E_{n,\chi,q} = \int_{X} [x]_q^n q^{-x} \chi (x) d\mu_{-q} (x)
.
$$

From this definition,  we derive

$$\split
\Cal E_{n,\chi,q} &= \int_{X} \chi (x)  [x]_q^n q^{-x}  d\mu_{-q}
(x)  \cr
 &= [d]_q^n \frac{[2]_q}{[2]_{q^d}} \sum_{a=0}^{d-1} \chi(a) (-1)^a \int_{\Bbb Z_p}
  [\frac{a}{d} +x]_{q^d} q^{-dx} d\mu_{-q^d}(x)\cr
  &= [d]_q^n \frac{[2]_q}{[2]_{q^d}} \sum_{a=0}^{d-1}   \chi(a) (-1)^a \Cal E_{n,q^d} (\frac{a}{d} ).
  \endsplit
$$

{\bf\centerline {\S 3. $q$-zeta function associated with $q$-Euler
numbers and polynomials}}

In this section we assume that $q\in\Bbb C$ with $|q|<1$. Let $F_q
(t,x) $ be the generating function of $\Cal E_{k,q} (x)$ as
follows:
$$
F_q (t,x ) = \sum_{n=0}^\infty \Cal E_{n,q}(x) \frac{t^n}{n!}.
$$
Then, we show that

$$\eqalignno{
F_q (t,x)  &=\sum_{m=0}^\infty \left( \frac{[2]_q}{(1-q)^m}
\sum_{l=0}^m \binom{m}{l} (-1)^l \frac{q^{lx}}{1+q^l}\right)
\frac{t^m}{m!} \cr
 &= [2]_q \sum_{m=0}^\infty \left( \frac{1}{(1-q)^m}
\sum_{l=0}^m \binom{m}{l} (-1)^l q^{lx}\sum_{k=0}^\infty (-1)^l
q^{kl} \right) \frac{t^m}{m!} &(22) \cr &= [2]_q \sum_{k=0}^\infty
(-1)^k \sum_{m=0}^\infty ( [k+x]_q^m \frac{t^m}{m!} ) =[2]_q
\sum_{k=0}^\infty (-1)^k e^{[k+x]_q t}.
  }
$$
Therefore, we obtain the following:

\proclaim{Theorem 8}  Let $F_q (t,x) $ be the generating function
of $\Cal E_{k,q} (x)$. Then we have
$$F_q (t,x) = [2]_q
\sum_{k=0}^\infty (-1)^k e^{[k+x]_q t} = \sum_{n=0}^\infty \Cal
E_{n,q} (x) \frac{t^n}{n!}.
$$
\endproclaim

From Theorem 8, we note that

$$ \Cal E_{k,q} (x)= \frac{d^k}{dt^k} F_q (t,x) |_{t=0} = [2]_q
\sum_{n=0}^\infty (-1)^n [n+x]_q^k.
$$

\proclaim{Corollary 9} For $k\in\Bbb Z_+$, we have
$$
\Cal E_{k,q} (x) = [2]_q \sum_{n=0}^\infty (-1)^n [n+x]_q^k  .$$
\endproclaim

\definition{Definition 10}
For $s\in\Bbb C$, we define $q$-zeta function as follows:
$$
\zeta_q (s,x) = [2]_q \sum_{n=0}^\infty \frac{(-1)^k}{[n+x]_q^s} .
$$
\enddefinition
Note that $ \zeta_q (-n , x) = \Cal E_{n,q} (x), $ for $n\in \Bbb
N \cup \{ 0\}$. Let $d$(=odd) be a positive integer. From the
generating function of $\Cal E_{n,q} (x)$, we derive
$$
\Cal E_{n,q} (x) =[d]_q^n  \frac{[2]_q}{[2]_{q^d}}
\sum_{a=0}^{d-1} (-1)^a \Cal E_{n, q^d} ( \frac{x+a}{d} ).
$$
Therefore we obtain the following:

\proclaim{Theorem 11} For $d$(=odd)$\in\Bbb N$, $n\in \Bbb Z_+$,
we have
$$
\Cal E_{n,q} (x) =[d]_q^n  \frac{[2]_q}{[2]_{q^d}}
\sum_{a=0}^{d-1} (-1)^a \Cal E_{n, q^d} ( \frac{x+a}{d} ).
$$
\endproclaim

Let $\chi$ be  the Dirichlet's character with conductor
$d$(=odd)$\in \Bbb N$ and let $F_{\chi,q} (t)$ be the generating
function of $\Cal E_{n,\chi,q}$ as follows:
$$F_{\chi,q} (t) =\sum_{n=0}^\infty \Cal E_{n,\chi, q}
\frac{t^n}{n!}.
$$
Then we see that
$$\split
F_{\chi,q} (t) &= \sum_{n=0}^\infty \left( [d]_q^n
\frac{[2]_q}{[2]_{q^d}} \sum_{a=0}^{d-1} \chi(a) (-1)^a  E_{n,
q^d} ( \frac{a}{d} ) \right) \frac{t^n}{n!}\cr
  &= \frac{[2]_q}{ [2]_{q^d} } \sum_{a=0}^{d-1} \chi(a)(-1)^a   \sum_{n=0}^\infty   E_{n,
q^d} ( \frac{a}{d} ) \frac{[d]_q^n t^n}{n!} \cr &=[2]_q
\sum_{a=0}^{d-1} \chi(a)(-1)^a \sum_{k=0}^\infty (-1)^k
e^{[kd+a]_q t}\cr &=[2]_q \sum_{n=0}^\infty \chi(n) (-1)^n
e^{[n]_q t}.
\endsplit
$$
Therefore, we obtain the following theorem:

\proclaim{Theorem 12} Let  $F_{\chi,q} (t)$ be the generating
function of $\Cal E_{n,\chi,q}$. Then we have
$$
F_{\chi,q} (t) = [2]_q \sum_{n=0}^\infty \chi(n) (-1)^n e^{[n]_q
t}=\sum_{n=0}^\infty \Cal E_{n,\chi,q} \frac{t^n}{n!} .
$$
\endproclaim

From Theorem 12, we derive
$$ \Cal E_{k,\chi,q}= \frac{d^k}{dt^k}
 F_{\chi,q}(t)|_{t=0}
=[2]_q \sum_{n=1}^\infty \chi(n) (-1)^n [n]_q^k .
$$

\proclaim{Corollary 13} For   $k\in\Bbb Z_+$, we  have
$$\Cal E_{k,\chi,q}
=[2]_q \sum_{n=1}^\infty \chi(n) (-1)^n [n]_q^k  .
$$
\endproclaim

For $s\in\Bbb C$, we define a $l$-series which interpolates the
modified  generalized $q$-Euler numbers attached to $\chi$ at a
negative integer as follows:
$$
l_q (s, \chi ) =[2]_q \sum_{n=1}^\infty \frac{\chi(n)
(-1)^n}{[n]_q^s}.
$$
From Corollary 13, we easily derive
$$
l_q (-n, \chi) =\Cal E_{n,\chi,q}, \ \ n \in\Bbb Z_+.
$$

\vskip 20pt

 Taekyun Kim

 EECS, Kyungpook National University,

  Taegu,  702-701, S. Korea

 e-mail:\text{tkim64$\@$hanmail.net, tkim$\@$knu.ac.kr}

\enddocument